\documentclass[10pt]{article}
\usepackage{graphics}
\usepackage{amssymb}
\usepackage{amsthm}

\newtheorem{theo}{\bf Theorem}

\def\qed{\hfill \rule{2.5mm}{2.5mm}}
\date{}
\begin{document}

\setcounter{footnote}{1}

\title{A Note on $(3,1)^*$-Choosable Toroidal Graphs
\thanks{This work is supported by NSFC of China,
 RFDP of Higher Education of China and Natural Sciences and
Engineering Research Council of Canada.}}
\author{Baogang Xu$^1$ and Qinglin Yu$^2$$^3$
\\ {\small $^1$ School of Mathematics and Computer Science}
\\ {\small Nanjing Normal University, Nanjing, China}
\\ {\small  Email:baogxu@pine.njnu.edu.cn}
\\ {\small  $^2$Center for Combinatorics, LPMC}
\\ {\small Nankai University, Tianjin, China}
\\ {\small  $^3$Department of Mathematics and Statistics}
\\ {\small Thompson Rivers University, Kamloops, BC, Canada}
\\ {\small  Email:yu@tru.ca}
}
\maketitle


\begin{abstract}
An $(L,d)^*$-{\it coloring} is a mapping $\phi$ that assigns a color
$\phi(v)\in L(v)$ to each vertex $v\in V(G)$ such that at most $d$
neighbors of $v$ receive colore $\phi(v)$.
A graph is called $(m,d)^*$-{\it choosable}, if $G$ admits an
$(L,d)^*$-coloring for every list assignment $L$ with $|L(v)|\geq
m$ for all $v\in V(G)$.
In this note, it is proved that every toroidal graph, which
contains no adjacent triangles and contains no 6-cycles and
$l$-cycles for some $l \in \{5,7\}$, is $(3,1)^*$-choosable.
\begin{flushleft}
\end{flushleft}
\end{abstract}

\section{Introduction}
\renewcommand{\baselinestretch}{2}
Graphs considered in this paper are finite, simple and undirected. A {\it toroidal} graph
$G=(V, E, F)$ is a graph embedded on the torus, where $V$, $E$ and
$F$ denote the set of vertices, edges and faces of $G$, respectively.

A face of an embedded graph is said to be {\it incident} with the
edges and vertices on its boundary. Two faces are {\it adjacent} if
they share a common edge. In particular, two adjacent 3-faces are
often referred as {\it adjacent triangles}. The {\it degree} of a
face $f$ of $G$, denoted by $d_G(f)$, is the number of edges
incident with it. Note that each cut-edge is counted twice in the
degree.  A $k$-{\it vertex} (or $k$-{\it face}) is a vertex (or a
face) of degree $k$, a $k^-$-vertex (or $k^-$-face) is a vertex (or
a face) of degree {\it at most} $k$, and a $k^+$-vertex (or
$k^+$-face) is defined similarly.  An $n$-face $f$ is called an
$(l_1,l_2,\ldots,l_n)$-face if the vertices incident with $f$ have
degree $l_1, l_2,\ldots,l_n$ sequentially. A cycle is called an
$m$-cycle if it is of length $m$. Undefined terms and notion follow
\cite{bb}.

    For each vertex $v\in V(G)$, we assign a set of
colors, $L(v)$ (called it {\it list}), to $v$. An $L$-{\it coloring} with
impropriety $d$ for non-negative integer $d$, or simply
$(L,d)^*$-{\it coloring}, is a mapping $\phi$ that assigns a color
$\phi(v)\in L(v)$ to each vertex $v\in V(G)$ such that at most $d$
neighbors of $v$ receive colore $\phi(v)$. For integers $m\geq
d\geq 0$, a graph is called $(m,d)^*$-{\it choosable}, if $G$ admits an
$(L,d)^*$-coloring for every list assignment $L$ with $|L(v)|\geq
m$ for all $v\in V(G)$. An $(m, 0)^*$-choosable graph is simply
called $m$-choosable.

    It is a hard problem to decide if a plane graph is 3-choosable,
even for triangle-free plane graphs. Thomassen proved that every
plane graph of girth at least 5 is 3-choosable \cite{10}. In
\cite{l3}, Voigt and Wirth constructed a family of triangle-free
plane graphs that is not 3-choosable. In \cite{hzbx}, it was
proved that every triangle-free plane graph containing no 8- and
9-cycles is 3-choosable. However, the 3-choosability of
triangle-free plane graph without 6- and 7-cycles is still open.

    The concept of list improper coloring was first introduced
    by \v{S}kreko\-vski \cite{ars}, and Eaton and Hull
\cite{aneth}, independently. They proved that every plane graph is
$(3,2)^*$-choosable and every outerplanar graph is
$(2,2)^*$-choosable.  \v{S}kreko\-vski \cite{ars2,ars1} investigated
the relationship between $(m, d)^*$-choosability and the girth in
plane graphs.  For instance, he proved every plane graph $G$ is $(3,
1)^*$-choosable if its girth, $g(G)$, is at least 4, and is $(2,
d)^*$-choosable if $g(G) \geq 5$ and $d\geq 4$. In \cite{aklzs}, it
was showed that every plane graph without 4-cycles and $l$-cycles
for some $l\in \{5,6,7\}$ is $(3,1)^*$-choosable.

    For toroidal graphs, Xu and Zhang \cite{bghz} proved that every
toroidal graph without adjacent triangles is $(4,1)^*$-choosable. In this note, we make the further restriction
$(3,1)^*$-choosability on toroidal graphs to improve Xu and Zhang's result to $(3,1)^*$-choosable.

    Let ${\cal G}$ denote the family of toroidal graphs containing
no adjacent triangles and containing no 6-cycles and $l$-cycles
for $l \in \{5,7\}$. The main result is to show that every graph
in ${\cal G}$ is $(3,1)^*$-choosable. In order to prove the main
theorem, we use the technique of ``discharging" to obtain several
forbidden configurations for the graphs in ${\cal G}$ and state as
a theorem below.

\renewcommand{\baselinestretch}{1.1}
\begin{theo}\label{th-1}
For every graph $G\in {\cal G}$, one of the following must hold:

$(1)$ $\delta(G)<3$.

$(2)$ $G$ contains two adjacent $3$-vertices.

$(3)$ $G$ contains a $(3, 4, 4)$-face.

$(4)$ $G$ contains a $(3, 4, 3, 4)$-face.
\end{theo}
\renewcommand{\baselinestretch}{2}

    As a consequence of the above result, we can prove the following

\renewcommand{\baselinestretch}{1.4}
\begin{theo}\label{th-2}
Every graph  in ${\cal G}$ is $(3,1)^*$-choosable.
\end{theo}
\renewcommand{\baselinestretch}{2}

\section{Proofs of the theorems}

    In the proof of Theorem \ref{th-1}, we use the technique of discharging.
In the beginning, each vertex $v$ is assigned a charge $\frac{d_G(v)}{3}-1$ and each face $f$ is
assigned a charge $\frac{d_G(f)}{6}-1$. By following the rules stated in the proof of the theorem, we will
redistribute the charges for the vertices and faces so that the new chages are nonnegative and the sume of the
new charges is still the same as before, which leads to a contradiction to Euler's formula.\\

\noindent \textbf{Proof of Theorem \ref{th-1}:} Assume to the contrary that
the theorem does not hold. Let $G$ be a connected toroidal graph in ${\cal G}$
satisfying $\delta(G)\geq 3$, every
3-vertex is adjacent to only $4^+$-vertices, and $G$ contains
neither $(3, 4, 4)$-faces nor $(3, 4, 3, 4)$-faces.

    Recall that we can rewrite Euler's formula $|V|+|F|-|E|=0$ for toroidal graphs as
\begin{equation}
\label{equ1} \sum_{v\in V(G)}\{{d_G(v)\over 3}-1\}+\sum_{f\in
F(G)}\{{d_G(f)\over 6}-1\}=0
\end{equation}

Defining a charge function $\omega$ on $V(G)\cup F(G)$ by letting
$\omega(v)={d_G(v)\over 3}-1$ if $v\in V(G)$ and
$\omega(f)={d_G(f)\over 6}-1$ if $f\in F(G)$. Then the total sum
of the charges, $\sum_{x\in V(G)\cup F(G)}\omega(x)$, is zero. \\

For two elements $x$ and $y$ of $V(G)\cup F(G)$, we use
$W(x\rightarrow y)$ to denote the charge transferred from $x$ to $y$.

\vskip 10pt

\noindent{\it Case 1.} $G\in {\cal G}$ contains neither 5- nor
6-cycles.

By the choice of $G$, it is easy to have the following observation.

$(O_{1,1})$ $G$ contains no 5- and 6-faces, no adjacent
$4^-$-faces.

Let $v$ be a $k$-vertex of $G$ and $f$ a $3$- or $4$-face incident
with $v$. Denote the number of $3$- or $4$-faces incident with $v$
by $r_v$. Then it is not hard to see $r_v \leq \lfloor{k \over
2}\rfloor$.

    The new charge function $\omega'(x)$ is obtained by
following discharging rules given below:

\vskip 5pt
\renewcommand{\baselinestretch}{0}
\textbf{$(R_{1,1})$} For all $v$ and $f$, $W(v\rightarrow
f)=\frac{1}{6}$ if $k=4$; $W(v\rightarrow f)=\frac{1}{3}$ if $k \geq
5$.

\textbf{$(R_{1,2})$} Each $7^+$-face transfers $\frac{1}{42}$ to each of
its adjacent $4^-$-faces.
\renewcommand{\baselinestretch}{2}

\vskip 5pt
    Reader is reminded that a face may be adjacent to another face {\it multiple}
    times.  Now, we ought to prove that $\omega'(x) \geq 0$ for any $x \in V(G)\cup F(G)$.

If $k=3$, then $\omega'(v)=\omega(v)=0$.

If $k=4$, then $r_v\leq 2$ and  $\omega'(v)\geq
\omega(v)-\frac{r_v}{6}=\frac{2-r_v}{6}\geq 0$.

If $k\geq 5$, then $\omega'(v)\geq
\omega(v)-\frac{r_v}{3}=\frac{k-3-r_v}{3}\geq 0$ (note that $r_v\leq
2$ if $k=5$).

\vskip 5pt

    Let $f$ be an $h$-face of $G$.  If $h \geq 7$, then
$\omega'(f)\geq \omega(f)-h\cdot\frac{1}{42}=\frac{6h-42}{42}\geq 0$.

If $h$ = 3, by ($O_{1,1}$), $f$ is adjacent to three $7^+$-faces.
 Since $G$ contains no adjacent 3-vertices and contains
no $(3,4,4)$-face, $f$ is either incident with a $5^+$-vertex and a
$4^+$-vertex or incident with three $4^+$-vertices. In the former
case, $f$ receives at least $\frac{1}{3}$ from the $5^+$-vertex and
receives at least $\frac{1}{6}$ from another $4^+$-vertex, and hence
\begin{equation}\label{eqa-3-1}
\omega'(f)\geq \frac{-1}{2}+{1\over 3}+{1\over 6}+3\cdot{1\over
42}>0.\end{equation}
 In the latter case, $f$ receives at least
${1\over 6}$ from each of the vertices incident with it, and hence
\begin{equation}\label{eqa-3-1-2}
\omega'(f)\geq {-1\over 2}+3\cdot {1\over 6}+3\cdot{1\over 42}>0.
\end{equation}

If $h$ = 4, then $f$ is adjacent to four $7^+$-faces. Since $G$
contains neither adjacent 3-vertices nor $(3,4,3,4)$-faces, $f$ is
incident with at least two $4^+$-vertices. Furthermore, if $f$ is
incident with a 3-vertex, then $f$ is either incident with a
$5^+$-vertex or incident with three $4^+$-vertices. Therefore,
\begin{equation}\label{eqa-4-1}
\omega'(f)\geq \omega(f)+{1\over 3}+{1\over 6}+4\cdot{1\over
42}>0.\end{equation}

    Thus,  $\omega'(x)\geq 0$ for each $x\in V(G) \cup
F(G)$. By (\ref{eqa-3-1}), (\ref{eqa-3-1-2}) and (\ref{eqa-4-1}),
$\omega'(f)>0$ if $d_G(f) = 3,4$. If $G$ contains no 3- and
4-faces, then $\omega'(f) = \omega(f) >0$ for any face $f$. Therefore,
 $0< \sum_{x\in V(G)\cup F(G)}\omega'(x)=\sum_{x\in
V(G)\cup F(G)}\omega(x)=0$, a contradiction.

\vskip 15pt

\noindent{\it Case 2.} $G\in {\cal G}$ contains neither 6- nor
7-cycles.

By the choice of $G$, we have the following observations.

$(O_{2,1})$ $G$ contains no 6- and 7-faces, no adjacent 3-faces, and
no adjacent

\hskip 33pt 4-faces.

$(O_{2,2})$  No 5-face is adjacent to 3- or 4-faces.

$(O_{2,3})$ Each 3-face is adjacent to at most one 4-face and each
4-face is

\hskip 33pt adjacent to at most one 3-face.

\vskip 10pt

Let $v$ be a $k$-vertex and $f$ an $l$-face incident with $v$. Denote
the numbers of $4^-$-faces and $5^-$-faces incident with $v$ by $r_1$ and $r_2$, respectively.
By $(O_{2,2})$ and $(O_{2,3})$,  we can see that $r_1\leq
\lfloor{2k\over 3}\rfloor$ and $3\lceil
\frac{r_1}{2}\rceil + r_2 \leq k+1$.

The discharging rules are as follows:

\vskip 5pt
\renewcommand{\baselinestretch}{0}
{\bf $(R_{2,1})$} For $k = 4$, $W(v\rightarrow f)={1\over 6}$ if $3\leq
l\leq 4$; $W(v\rightarrow f)={1\over 18}$ if $l=5$.

{\bf $(R_{2,2})$} For $k \geq 5$, $W(v\rightarrow f)={1\over 4}$ if $l=3$;
$W(v\rightarrow f)={1\over 6}$ if $l=4$;

\hskip 33pt $W(v\rightarrow f)={1\over 18}$ if $l=5$.

{\bf $(R_{2,3})$} An $8^+$-face transfers ${1\over 24}$ to each of
its adjacent $5^-$-faces.

\vskip 8pt

    We now verify that $\omega'(x) \geq 0$ for any $x \in V(G)\cup F(G)$.

If $k = 3$, then $\omega'(v)=\omega(v)=0$.

If $k = 4$, then $r_1\leq 2$.  From $(O_{2,2})$, $r_1=2$ implies $r_2=0$
and $r_2\geq 2$ implies $r_1=0$. Hence $3r_1+r_2\leq 6$ and $\omega'(v)\geq \omega(v)-{r_1\over
6}-\frac{r_2}{18}=\frac{1}{3}-\frac{1}{18}(3r_1+r_2)\geq 0$.

If $k = 5$, then $r_1\leq 3$. If $r_1=3$, then $r_2=0$ and
$v$ is incident with a $4$-face. Hence
$\omega'(v)\geq \omega(v)-2\cdot{1\over 4}-{1\over 6}=0$. If
$r_1=2$, then $r_2\leq 1$ and thus  $\omega'(v)\geq
\omega(v)-2\cdot{1\over 4}-{1\over 18}>0$. If $r_1\leq 1$, then
$r_2\leq 2$ and thus $\omega'(v)\geq \omega(v)-{1\over 4}-{2\over 18}>0$.

If $k\geq 6$, since $G$ contains no adjacent triangles,
$v$ is incident with at most $\lfloor {k\over 2}\rfloor$ 3-faces.
If $r_1\leq  \lfloor {k\over 2}\rfloor$, then

$\begin{array}{rlllrrr}
    \omega'(v)&\geq  & \omega(v)-{r_1\over 4}-{r_2\over 18}={k-3\over 3}-{9r_1+2r_2\over
36}\\
              & = &{k-3\over 3}-{3(3\lceil{r_1\over 2}\rceil+r_2)+9\lfloor{r_1\over
2}\rfloor-r_2\over 36}\\
              &\geq &{12k-36-3(k+1)-9\lfloor{r_1\over 2}\rfloor+r_2\over 36}\geq
{9k-33-9\lfloor{k\over 4}\rfloor\over 36}\geq 0.
\end{array}$

\noindent If $r_1>\lfloor {k\over 2}\rfloor$, then $v$ is incident
with at least $r_1-\lfloor {k\over 2}\rfloor$ 4-faces and thus

$\begin{array}{rlllrrr} \omega'(v) & \geq &\omega(v)-{1\over
4}\lfloor {k\over 2}\rfloor-{1\over 6}(r_1-\lfloor {k\over
2}\rfloor)-{r_2\over 18}\\
    & = & {k-3\over 3}-{9\lfloor {k\over
2}\rfloor+6r_1-6\lfloor {k\over 2}\rfloor+2r_2\over 36}\\
        & = & {k-3\over 3}-{3\lfloor {k\over 2}\rfloor+2(3\lceil{r_1\over
2}\rceil+r_2)+6\lfloor{r_1\over 2}\rfloor\over 36} \\
        & \geq & {12k-36-3\lfloor {k\over 2}\rfloor-2(k+1)-6\lfloor{r_1\over
2}\rfloor\over 36}\geq {10k-38-3\lfloor {k\over
2}\rfloor-6\lfloor{k\over 3}\rfloor\over 36}\geq 0.
\end{array}$

\vskip 5pt

    Let $f$ be an $h$-face of $G$.

    If $h\geq 8$, then, by $(R_{2,3})$,
$\omega'(f)\geq \omega(f)-{h\over 24}={3h-24\over 24}\geq 0$.

    If $h = 5$, then  $f$ is incident with at least three $4^+$-vertices
(note that $G$ contains no adjacent 3-vertices) and,  by $(R_{2,1})$
and $(R_{2,2})$,  each of these $4^+$-vertices transfers at least
${1\over 18}$ to $f$ and hence $\omega'(f) \geq \omega(f)+{3\over
18}=0$ if $f$ is not adjacent to $8^+$-faces and
\begin{equation}\label{eqa-5-2} \omega'(f)=\omega(f)+{3\over
18}+{1\over 24}>0 \ \ \ \ \mbox{if} \ f \mbox{ is adjacent to at
least one} \ \  \mbox{$8^+$-face}.\end{equation}

If $h$ = 4, $f$ is incident with at least two $4^+$-vertices and is
adjacent to at least three $8^+$-faces.  Thus
\begin{equation}\label{eqa-4-2}\omega'(f)\geq \omega(f)+2\cdot {1\over 6}+3\cdot {1\over
24}>0.\end{equation}

If $h$ = 3, $f$ is adjacent to at least two $8^+$-faces of which
each transfers ${1\over 24}$ to $f$. Therefore, $\omega'(f)\geq
{-1\over 2}+{1\over 4}+{1\over 6}+2\cdot{1\over 24}=0$ if
$f$ is incident with a $5^+$-vertex and another $4^+$-vertex, and
$\omega'(f)\geq {-1\over 2}+3\cdot {1\over 6}+2\cdot{1\over 24}>0$
if $f$ is incident with three $4^+$-vertices.

\vskip 5pt

    By (\ref{eqa-4-2}), $\omega'(f)>0$ for
every 4-face $f$. From (\ref{eqa-5-2}), $\omega'(f)>0$ for every
5-face $f$ adjacent to some $8^+$-faces. If $G$ contains no
4-faces, then every 3-face is adjacent to three $8^+$-faces that
yields $\omega'(f)>0$ for any 3-face $f$. If $G$ contains no
$5^-$-face adjacent to $8^+$-faces, then $\omega'(f)>0$ for any
$8^+$-face $f$. So, $0<\sum_{x\in V(G)\cup F(G)}\omega'(x)=\sum_{x\in
V(G)\cup F(G)}\omega(x)=0$. This contradiction leads to the proof of Case 2
and thus Theorem \ref{th-1}.  \qed

\vskip 15pt

\noindent{\bf Proof of Theorem \ref{th-2}:} Assume to the contrary.
Let $G$ be a counterexample with the fewest vertices, i.e., there
exists a list assignment $L$ with $\mid$$L(v)$$\mid$=3 for all $v\in
V(G)$ such that $G$ is not $(L,1)^*$-choosable, but any proper
subgraph of $G$ is.

  If $\delta(G)<3$, let $v$ be a $2$-vertex of $G$. Then, $G-v$ is
$(3,1)^*$-choosable by the choice of $G$. Since in any
$(L,1)^*$-coloring of $G-v$, there must exist a color in $L(v)$ that
is not used by any neighbors of $v$, any $(L,1)^*$-coloring of $G-v$
can be extended to a $(L,1)^*$-coloring of $G$, a contradiction. So
we assume that $\delta(G)\geq 3$.

If $G$ contains two adjacent 3-vertices, say $u$ and $v$, then by
the choice of $G$, $G-\{u,v\}$ is $(3,1)^*$-choosable. In any
$(L,1)^*$-coloring of $G-\{u,v\}$, there exists a color in $L(u)$
that is not used by any neighbors of $u$ in $G-\{u,v\}$, and the
same holds for $v$. Applying the same argument as the above, we see that $G$ is
$(L,1)^*$-choosable, a contradiction.

    Suppose that $G$ contains a $(3, 4, 4)$-face $f$ with the boundary $xyzx$,
say, $d_G(x)=3$ and $d_G(y)=d_G(z)=4$. Let
$H=G-\{x,y,z\}$. By the choice of $G$, $H$ admits an
$(L,1)^*$-coloring $\phi$. For $w\in \{x,y,z\}$, let
$L'(w)=L(w)\setminus\{\phi(u)|u\in N_H(w)\}$. Then, $|L'(x)|\geq
2$, $|L'(y)|\geq 1$ and $|L'(z)|\geq 1$. If $L'(y)=L'(z)$, then
color $y$ and $z$ with a same color $\gamma$ in $L'(y)$ and color
$x$ with a color in $L'(x)\setminus \{\gamma\}$. If $L'(y)\neq
L'(z)$, then color $y$ with a color $\alpha\in L'(y)\setminus
L'(z)$, color $z$ with a color in $L'(z)$, and color $x$ with an
arbitrary color in $L'(x)$. In either case, we obtain an
$(L,1)^*$-coloring of $G$, a contradiction.

By Theorem \ref{th-1}, we may assume that $G$ contains a
$(3, 4, 3,4)$-face $f$ with the boundary $xyzux$. By symmetry, we
assume that $d_G(x)=d_G(z)=3$ and $d_G(y)=d_G(u)=4$. Let
$F=G-\{x,y,z,u\}$. By the choice of $G$, $F$ admits an
$(L,1)^*$-coloring $\psi$. For $w\in \{x,y,z,u\}$, let
$L'(w)=L(w)\setminus\{\psi(v)|v\in N_F(w)\}$. Then, $|L'(x)|\geq
2$, $|L'(z)|\geq 2$, $|L'(y)|\geq 1$ and $|L'(z)|\geq 1$. It is
easy to verify that $xyzu$ admits an $(L',1)^*$-coloring. This
together with $\psi$ yields an $(L,1)^*$-coloring of $G$. This
contradiction completes the proof of Theorem \ref{th-2}. \qed

\vskip 20pt

\noindent \title{\Large\bf Acknowledgments} \maketitle
 The authors are indebted to the anonymous referees for their constructive suggestions.

\vskip 20pt

\renewcommand{\baselinestretch}{0}

\end{document}